\input amstex
\documentstyle{amsppt}

\def\e{\epsilon}
\def\ZZ{\Bbb Z}
\def\NN{\Bbb N}

\def\PP{\Bbb P}

\def\Dm,x{{D_{m,x_{1},\ldots,x_{18}}}}


\topmatter
\title On Genera of Smooth Curves in Higher Dimensional Varieties \endtitle
\author Jungkai Alfred Chen \endauthor
\address{ Jungkai Alfred Chen \hfill {} \linebreak
\hglue 13.6pt Department of Mathematics \hfill {} \linebreak
\hglue 13.6pt University of California \hfill {} \linebreak
\hglue 13.6pt Los Angeles, CA 90095-1555}  \endaddress
\email  jachen\@math.ucla.edu \endemail
\abstract
We prove that for any smooth projective variety $X$ of dimension $\geq 3$, 
there exists an integer $g_0=g_0(X)$, such that for any integer $g \geq g_0$, 
there exists a smooth curve $C$ in $X$ with $g(C) = g$.
\endabstract
\endtopmatter

\document
\heading Introduction \endheading
It's a very elementary fact that for any integer $g \geq 0$, 
there exists a smooth curve 
$C \subset \PP^3$ of genus $g$.
It is interesting to ask what the analogous situation is 
when projective space is replaced by an arbitary
smooth projective variety $X$ of dimension $\geq 3$. 
It can happen of course that a given variety $X$ 
contains no curves of small genus.
For example, abelian varieties contain no rational curves,
and Clemens \cite{C} has shown that 
on a generic hypersurface of degree $d \text{ in } \PP^n$,
there are no curves of genus $\leq \frac{d-2n+1}{2}$.
So the natural question is 
whether all sufficiently large genera are realized by smooth curves. 
Our main result states that this is indeed the case.

\proclaim{Theorem 1}
Let $X$ be a smooth projective variety of dimension $n \geq 3$.
Then there exists an integer $g_0=g_0(X)$ such that 
for any integer $g \geq g_0$, 
there exists a smooth curve $C \subset X$ of genus $g$.
\endproclaim

We also prove an analogous statement for the geometric genus of 
nodal curves on surfaces.

The paper is organized as follows. 
In section 1, 
we prove a lemma that a certain type of function will represent
all large enough integers. 
In section 2, 
we construct a smooth surface $S$ in 
any smooth projective variety of dim $\geq 3$ such that 
$Pic(S)$ is very large.
We choose divisors in $S$ which give us various genera. 
In section 3, 
we prove the analogous statement for surfaces that 
any large enough integer
can be realized as the geometric genus of some nodal curves on a given surface. 

We'd like to thank Robert Lazarsfeld for his encouragement and suggestion.
We are indebted to Joe Harris and Claire Voisin for pointing out 
the importance of surfaces with large Picard groups. 
We would also thank Richard Elman, Mark Green and Christopher Hacon 
for helpful discussion.

\heading 1. A Numerical Lemma  \endheading
\proclaim{Lemma 1}
Let $$f(x_{1},\ldots,x_{18}) = \sum_{i=1}^{9} ax_{i}^{2} +
 \sum_{i=10}^{18} bx_{i}^{2} + \sum_{i=1}^{18} cx_{i},$$
where $a,b \in \NN \text{, } c \in \ZZ \text{, and } |a-b| = 1 \text{ or } 2$.
Then there exists an integer $m_{0}=m_{0}(a,b)$ 
 such that for any even integer $m \geq m_{0}$,
 $f(x_{1},\ldots,x_{18}) = m$ has an integral solution.
\endproclaim

\demo{Proof}

1. suppose $(a,b)=1$: 

There exists an $n_{0} \in \ZZ$  such that for any integer $n \geq n_{0}$, 
$n=ar +bs$  for some $r,s \in \NN$.
Furthermore, every positive integer is sum of four squares. So
$$  n=a \sum_{i=1}^{4} \bar{x}_{i}^{2}+b\sum_{i=5}^{8} \bar{x}_{i}^{2} $$ 
for some $\bar{x}_{1},\ldots,\bar{x}_{8} \in \ZZ$.
Hence for any even integer $m \geq 2n_{0}$
$$ m=f(\bar{x}_{1},\ldots,\bar{x}_{4},-\bar{x}_{1},\ldots,-\bar{x}_{4},0,
\bar{x}_{5},\ldots,\bar{x}_{8},-\bar{x}_{5},\ldots,-\bar{x}_{8},0) $$ 
for some $\bar{x}_{1},\ldots,\bar{x}_{8} \in \ZZ$. 

2. suppose $(a,b)=2$:

(i) suppose  $m \equiv 0 \pmod{4}$: 

We write $\frac{m}{4}$ as a combination of $\frac{a}{2}$ and 
$\frac{b}{2}$ for $m \gg 0$ . So similiarly, 
$$ m=f(\bar{x}_{1},\ldots,\bar{x}_{4},-\bar{x}_{1},\ldots,-\bar{x}_{4},0,
\bar{x}_{5},\ldots,\bar{x}_{8},-\bar{x}_{5},\ldots,-\bar{x}_{8},0) $$ 
for some $\bar{x}_{1},\ldots,\bar{x}_8 \in \ZZ$. 

(ii) suppose  $m \equiv 2 \pmod{4}$:
 
Then we have $a+b \equiv 2 \pmod{4}$ and $m-a-b \equiv 0 \pmod{4}$ 
Hence for $m \gg 0$ 
$$ m-a-b = f(\bar{x}_{1},\ldots,\bar{x}_{4},-\bar{x}_{1},\ldots,-\bar{x}_{4},0,
	 \bar{x}_{5},\ldots,\bar{x}_{8},-\bar{x}_{5},\ldots,-\bar{x}_{8},0).$$
$$m  = f(\bar{x}_{1},\ldots,\bar{x}_{4},-\bar{x}_{1},\ldots,-\bar{x}_{4},1, 
    \bar{x}_{5},\ldots,\bar{x}_{8},-\bar{x}_{5},\ldots,-\bar{x}_{8},-1).$$
for some $\bar{x}_{1},\ldots,\bar{x}_8 \in \ZZ$. 
This proves the lemma. 
\hfill \qed 
\enddemo

\remark{Remark 1} 
From the above proof, $\sum_{i=1}^{4} a\bar{x}_i^2 \leq m \text{ and }
                       \sum_{i=5}^{8} b\bar{x}_i^2 \leq m$.
So $\bar{x}_i \leq \sqrt{m}$.
We can modify Lemma 1 as:
 
Given f, $\forall m \gg 0$, 
$m$ can be represented by $f(x_1,\ldots,x_{18}) \text{ with }
|x_i| \leq \sqrt{m}$ for all $i$. 
\endremark

\heading 2. Main Theorem  \endheading
It's enough to prove the theorem when dim $X = 3$.
So let $X$ be a smooth projective threefold. 
We first construct curves $C_{1} \text{ and } C'_{1} \subset X$ 
($C_{1}$ smooth but $C'_{1}$ singular) 
whose arithmetic genera differ by one.
Let $H_{i} (i=1 \ldots 5)$ denote 
linearly equivalent very ample divisors in $X$.
Consider $$C_{1}=(H_{1} \cap H_{2}) \cup (H_{3} \cap H_{4}) \text{, }
C'_{1}=(H_{1} \cap H_{2}) \cup (H_{3} \cap H_{5}).$$
Choose $H_{i}$'s properly so that $C_{1}$ is smooth but disconnected
(with $H_{1} \cap H_{2} \text{ and } H_{3} \cap H_{4}$ as its two components),
and $C'_{1}$ is connected and has a simple node at 
$H_{1} \cap H_{2} \cap H_{3} \cap H_{5}$. 
Such $C_{1}$ and $C'_{1}$ are numerically equivalent in $X$, 
and their arithmetic genera differ by one.
By moving the very ample divisors, 
we can construct more curves,
 $C_{1}, \ldots, C_{18}$, such that
 
1. They all are disjoint. 

2. For all divisors $D \subset X$, $\forall$ $1\leq i,j \leq 18$,
        $$C_i \cdot D=C_j \cdot D.$$ 

3. $p_{a}(C_{i})=p_{a}(C_{j})+1$,
        $\forall$ $1 \leq i \leq 9 \text{, } 10 \leq j \leq 18$.

Fix next a very ample divisor $H$ in $X$ 
and a smooth surface $S \in |nH|, (n \gg 0)$, 
 containing all curves constructed above.
Computing on $S$ for any of our curves $C$:
$$2p_{a}(C)-2  =  K_{X} \cdot C+nH \cdot C + C \cdot_{S} C.$$
Let $$C_{i} \cdot C_{i}=-a \text{, for all } 1 \leq i \leq 9,$$
then $$C_{i} \cdot C_{i}=-(a+2), \text{ for all } 10 \leq i \leq 18.$$
Since we choose $S \in |nH|$ with $n \gg 0$, 
We may assume that $a > 0$.

Next, let $H_S$ denote the restriction of 
the very ample divisor $H$ in $X$ to $S$. 
Let
$$D_{m,x_{1},\ldots,x_{18}}  \equiv  K_S+ mH_S+\sum_{i=1}^{18} x_{i}C_{i}, $$
$$D_{m}  \equiv  K_S+ mH_S. $$

The plan is this: 
We show that for $|x_{i}|$ small compared to $m$, 
the linear series $| D_{m,x_{1},\ldots,x_{18}} |$ is very ample.
By using lemma 1, 
we show that by choosing suitable $x_{i}$ 
we can find curves in their linear series 
of all possible genera sufficiently close to $p_{a}(D_{m})$.
Then we let $m$ vary.

\proclaim{Lemma 2}
For any $e>0$,
there exists positive integer $m_1=m_1(e)$ such that  
$|\Dm,x |$ is very ample 
$\forall (m \geq m_{1}, |x_{i}| \leq e \sqrt{m})$.
In particular, there exist smooth curves in $|\Dm,x |$.
\endproclaim

\demo{Proof} 
We prove the very ampleness of adjoint linear series 
by Reider's Theorem \cite{L}.
We need to show that for any irreducible curve $\Gamma \subset S,$
$$( \Dm,x -K_{S}) \cdot \Gamma \geq 3$$
and
$$( \Dm,x-K_{S})^2 \geq 10.$$ 
We can find an integer $t$ such that
$| tH_S -\sum_{i=1}^{18}\e_iC_i|$
are very ample for any $\e_i$=1, 0, or -1,
and hence $| tH_S -\sum_{i=1}^{18}\e_iC_i|$ 
has positive intersection with any effective curve.
If we pick $m$ such that $m > 18te\sqrt{m} \geq 18t |x_i|$,
easy computation shows that these are indeed the case.
\hfill \qed
\enddemo

\demo{Proof of Main Theorem} 
Let 
$$ g_m=p_a(D_m)$$ 
$$ S_m  = \{  p_a(\Dm,x)|m,x_i \in \ZZ \} =
          \{ g_m -\frac{1}{2} f(x_1,\ldots,x_{18}) \}$$ 
where $$f(x_1,\ldots,x_{18})=\sum_{i=1}^9(a+2)x_{i}^2+
\sum_{i=10}^{18}ax_{i}^{2}+\sum_{i=1}^{18}m(H \cdot C_i)x_i$$
satisfies the condition of lemma 1.

We first show that for any integer $z \leq g_{m-1}$ , 
$z$ will be in $ S_m.$
This can be done by Lemma 1 and letting $m \gg 0$
because $$ g_m -z \geq g_m-g_{m-1}=mH_S^2+o(1).$$

Secondly , we want to show that 
for any integer $z, g_{m-2} \leq z \leq g_{m-1}$,
$z$ is actually
genus of some smooth curves if $m \gg 0$.
Choose $e$ so that 
$e^{2}m \geq 2(g_{m}-g_{m-2}),\forall m \gg 0$.
Hence for any integer $z, g_{m-2} \leq z \leq g_{m-1}$,
$z=p_a(\Dm,x)$ for some $|x_i| \leq \sqrt{2(g_m-z)} \leq e\sqrt{m}$~[Remark 1].
By Lemma 2, such $z$ can be realized as genus of some smooth curves.

Finally, let $m$ vary.
Every integer $g $ lies inside some interval
$[g_{m-2}, g_{m-1}].$ 
This completes the proof.
\hfill \qed
\enddemo

\heading 3. Variants, and Example \endheading
We can't expect to realize all large genera 
by smooth curves on a general surface. 
The best thing we can hope is allowing curves to have some nodes. 
\proclaim{Theorem 2}
Let $S$ be a smooth surface. 
Then there exists an integer $g_{0}=g_{0}(S)$ such that 
for any integer $g \geq g_{0}$, 
there exists nodal curves with geometric genus $=g$ 
\endproclaim

\demo{Proof} 
Let $H$ be a very ample divisor on $S$.
Let $$g_m = p_a(K_S+mH),$$
    $$d_m = g_m - g_{m-1}.$$ 
It is enough to show that for $m \gg 0$, 
and $\forall r \text{, } 1 \leq r \leq d_m$, 
there exists a nodal curve $C \in |K_X+mH|$ with exactly $r$ nodes.

Let $C_1 \in |[\frac{m}{4}] H|$ be a smooth curve in $X$.
Let $\Sigma= \{x_1,\ldots,x_r\}$ be $r$ distinct points in $C_1$.
Now consider the blowup of $X$ along $\Sigma$
$$f: \tilde{X}=Bl_{\Sigma}(X) \rightarrow X.$$
Let $E_i$'s denote the exceptional divisors over $x_i$'s. 
Let $$L_{m,r} = mf^{*}H-3(E_1 + \ldots + E_r).$$
Then $$K_{\tilde{X}} + L_{m,r} \equiv 
   f^{*}( K_X+mH ) -2(E_1 + \ldots + E_r).$$  
$$L_{m,r} \equiv (m-3[\frac{m}{4}])f^*H + 3 \tilde{C_1},$$ 
where $\tilde{C_1}$ denotes the proper transformation of $C_1$.

Computation shows that, for $m$ large enough, 
$L_{m,r}$ has both self-intersection and intersection with effective
curves positive enough. 
So $|K_{\tilde{X}} + L_{m,r}|$ is very ample by Reider's theorem. 
Then we can pick an curve $\tilde{C} \in |K_{\tilde{X}} + L_{m,r}|$ 
such that
$\tilde{C}$ is smooth and $\tilde{C}$ intersects $E_i$ at 2 distinct points 
$\forall$ $1 \leq i \leq r.$
Hence $C:=f_{*}(\tilde{C})$ is a nodal curve in $|K_X+mH|$ with $r$ nodes 
at $x_1,ldots,x_r$.
This completes the proof.
We left the details to the readers.
\hfill \qed
\enddemo

Turning back to higher dimensional case,
we conclude with an elementary example
to show that the least integer $g_0$ appearing in Theorem 1 can be
arbitary large:

\example{Example}
 Let $X=C_{1}\times C_{2} \times C_{3}$, 
where $C_i$'s are smooth projective curves.
There are natural projection $p_i$'s to each $C_i$'s.
For any smooth irreducible ecurve $C$ in $X$,
not all the projections of $C$ are points.
Hence we have:
$$ g(C) \geq \min(g(C_{1}),g(C_{2}),g(C_{3})) $$ 
Hence there is no uniform $g_{0}$
for all projective varieties of same dimension. 
\endexample
\bigskip


\enddocument
\end